\documentclass[reqno,fleqn]{amsart}
\usepackage{amsfonts,amsthm,amssymb}
\usepackage{lipsum}
\usepackage[utf8]{inputenc}
\usepackage[english]{babel}

\usepackage{lmodern,textcomp}
\usepackage{amssymb}
\usepackage{amsthm}
\usepackage{amsmath}
\usepackage{mathtools}
\usepackage{latexsym}
\usepackage{tikz}
\usepackage{fancyvrb}
\usepackage{epsfig}
\usepackage{graphicx}
\usepackage{amsmath}
\usepackage{latexsym}
\usepackage{tikz}
\usepackage{rotating}
\usetikzlibrary{positioning}
\usepackage{subcaption}
\usepackage{datetime}
\newtheorem{theorem}{Theorem}[section]

\newtheorem{cor}[theorem]{Corollary}
\newtheorem{definition}[theorem]{Definition}
\newtheorem{conj}{Conjecture}
\newtheorem{exm}{Example}

\newtheorem{rem}[theorem]{Remark}

\newtheorem{note}[theorem]{Note}

\title[Finding Non-Distance Magic Graphs using nbh chains]{Finding Non-Distance Magic Graphs using neighbourhood chains}
\author{\sc V. Vilfred Kamalappan and Sajidha P} 
\address{Department of Mathematics, ~Central ~University ~of~ Kerala, Kasaragod, ~India.}

\email{vilfredkamalv@cukerala.ac.in}
\email{sajinada555@gamil.com}


\date{}

\begin{document} 
	
\begin{abstract} Let $G$ be a graph of order $n$ and $N = \{N(u_{i})\}^k_{i=1}$ be a sequence of neighbourhood(nbh)s in $G$ where $N(u)$ = $\{v\in V(G):$ $uv\in E(G)\}$. \emph{Nbh sequence graph $H$ of} $N$ in $G$ is defined as the union of all induced subgraphs of closed nbh $N[u_{i}]$ in $G$, $1 \leq i \leq k$, $k\in\mathbb{N}$. A labeling $f: V(G) \rightarrow \left\{1,2,\ldots,n\right\} $ is called a \emph{Distance Magic Labeling (DML)} of $G$ if ~ ${\sum_{v \in N(u)}} f(v) $ is a constant for every $u\in V(G)$. $G$ is called a \emph{Distance Magic graph (DMG)} if it has a DML, otherwise it is called a \emph{Non-Distance Magic (NDM)} graph. In this paper, we define nbh walk, nbh trial, nbh path or nbh chain, nbh cycle, nbh sequence graph and nbh chains of Type-1 (NC-T1) and Type-2 (NC-T2). NC-T2 is formed on two NC-T1 of same length. We prove that (i) for $k \geq 2$ and $n \geq 3$, cylindrical grid graph $P_{k} \Box C_{n}$ contains NC-T2, $k,n \in \mathbb{N}$; (ii) graph containing NC-T1  of even length is NDM and (iii) partially settle a conjecture that graphs $P_m \Box C_n$ are NDM when $n$ is even, $m \geq 2$, $n \geq 3$ and $m,n\in\mathbb{N}$.    
\end{abstract}
	\maketitle
	
\section{Introduction}
	
	In 1987, Vilfred \cite{v87} defined \emph{Sigma labeling, Sigma partition} and \emph{Sigma labeled graphs} and develpoed its theory \cite{v96}. The same labeling was independently defined in 2003 \cite{mrs} as \emph{1-distance magic vertex labeling} and the term `Distance Magic Labeling' was used in a 2009 article \cite{sf}. Hereafter, we use the term `Distance Magic Labeling' since the motivation to define Sigma labeling of graphs is the construction of Magic Squares (See pages 1, 97-99 of \cite{v96}.).
	
	 In 1996, while studying sigma labeling of graphs, Vilfred proved that cylindrical grid graphs $P_2 \Box C_n$ for $n \geq 3$, and $P_k \Box  C_n$ for $n$ = 3,4 and $k \geq 2$ are Non-Distance Magic (NDM) and proposed the following  conjecture in \cite{v96} (pages 11-13 in \cite{v96}).
	\begin{conj}  \label{ca1} {\rm \cite{v96} For $m \geq 2$, $n \geq 3$ and $m,n\in\mathbb{N}$, cylindrical grid graphs $P_m \Box C_n$ are Non-Distance Magic (NDM). } 
	\end{conj}
	
	In this paper, we define neighbourhood(nbh) walk, nbh trial, nbh path or nbh chain, nbh cycle similar to the definition of walk, trial, path, cycle in graphs and also nbh sequence graphs. These are used to develope the theory in this paper. We also define nbh chains of Type-1 (NC-T1) and Type-2 (NC-T2). NC-T2 is formed on two NC-T1 of same length. We prove that (i) for $k \geq 2$ and $n \geq 3$, cylindrical grid graph $P_{k} \Box C_{n}$ contains NC-T2, $n,k \in \mathbb{N}$; (ii) graph containing NC-T1  of even length is NDM and (iii) for $m \geq 2$ and $n \geq 3$, $P_{2m} \Box C_n$ are NDM and thereby, we partially settle Conjecture \ref{ca1}, $m,n\in\mathbb{N}$.
	
	Effort to settle Conjecture \ref{ca1} is the motivation to do this research work. 
	
	Through out this paper, we consider only finite undirected simple graphs and for all basic ideas in graph theory, we follow \cite{dw01}. We present here a few definitions and results which are needed in the subsequent sections.
		
	\begin{definition} \cite{v87} \quad	Let $G$ be a simple graph of order $n$. Then a  bijective mapping $f: V(G) \rightarrow \left\{1,2,\ldots,n\right\} $ is called a \emph{Distance Magic Labeling (DML)} or \emph{$\Sigma$-labeling (Sigma labeling) of $G$} if ~ ${\sum_{v \in N(u)}} f(v) $ = $S$ is a constant for all  $u\in V(G)$ where $N(u)$ = $\{v\in V(G): uv\in E(G)\}$. Graph $G$ is called a \emph{Distance Magic graph (DMG)} if it has a DML, otherwise it is called a \emph{Non-Distance Magic (NDM) graph} and we call $S$ as \emph{Distance Magic constant}.
	\end{definition}	  
	
	\begin{definition}\cite{dw01} A \emph{path} is a simple graph whose vertices can be ordered so that two vertices are adjacent if and only if they are consecutive in the list. The (unlabeled) path with $n$ vertices is denoted by $P_n$.
	\end{definition}
	
	\begin{definition}\cite{dw01} 
		A closed path on $n$ vertices is called a \emph{cycle} and is denoted by $C_{n}$, $n \geq 3$ and $n\in\mathbb{N}$.
	\end{definition}
	
	\begin{definition}\cite{dw01} The \emph{cartesian product of graphs} $G$ and $H$, written $G \Box H$, is the
		graph with vertex set $V(G)$ $\times$ $V(H)$ specified by putting $(u, v)$ adjacent to
		$(u', v')$ if and only if (i) $u = u'$ and $vv'\in E(H)$, or (ii) $v = v'$ and $uu'\in E(G)$.
	\end{definition}
	
	\begin{definition}\cite{dw01} The \emph{union of two disjoint graphs} $G$ and $H$, written $G \cup H$, is
		the graph with vertex set $ V(G) \cup V(H)$ and edge set $E(G) \cup E(H)$.
	\end{definition}
	
\begin{theorem}  \label{a0} {\rm {[Theorem 2.6 in \cite{v96}]}\quad For $m \geq 2$, graphs $ P_m \Box C_3 $  and  $ P_m \Box C_4 $ are NDM, $ m \in \mathbb{N}$. \hfill $\Box$}
\end{theorem}
	
	In \cite{mrs}, the following condition is given to identify NDM graphs.
	
	\begin{theorem}\cite{mrs} \label{a2} {\rm If a simple graph $G$ contains two distinct vertices $u$ and $v$ $\ni$ $ |N (u) \cap N (v)|$ = $deg(v) - 1$ = $deg(u) - 1 $, then $G$ is NDM. \hfill $\Box$ }
	\end{theorem}
		
	While trying to prove Conjecture  \ref{ca1}, we could notice that for $m \geq 2$, $n \geq 3$ and $ m,n \in \mathbb{N}$, the cyclindrical grid graph $P_{m} \Box C_{n}$ satisfies a particular condition which is a more general case of condition given in Theorem \ref{a2} and is used here as a tool to identify a large families of graphs as NDM graphs. 
		
	\section{nbh walk, nbh trial, nbh chain, nbh cycle, nbh sequence graphs}
	
	Here, we define neighbourhood(nbh) walk, nbh trial, nbh path or nbh chain, nbh cycle similar to the definition of walk, trial, path, cycle in graphs and also nbh sequence graphs. These are used to develope the theory in this paper. 
		
	A \emph{nbh sequence} in a graph $G$ is a sequence of nbhs in $G$. That is a nbh squence is a sequence in which each element/term is a nbh of a vertex in $G$.
	
	\begin{rem} \quad  In general, in any nbh sequence $\{N(u_{i})\}^k_{i=1}$ in a graph $G$, we consider $k \geq 2$ and $|N[u_{i}]| > 1$, $\forall$ $i$, unless stated otherwise, $u_{i}\in V(G)$, $1 \leq i \leq k$ and $i\in\mathbb{N}$. 
	\end{rem}

	\begin{definition} \quad  A nbh sequence $\{N(u_{i})\}^k_{i=1}$ in a graph $G$ is called a  \emph{nbh walk} if consecutive terms(nbhs) are having common point(s), $u_i\in V(G)$, $1 \leq i \leq k$. We denote the nbh walk $\{N(u_{i})\}^k_{i=1}$ as $N_1 N_2 \cdots  N_k$ where $N_i$ = $N(u_i)$, $u_i\in V(G)$, $k \geq 2$ and $i\in\mathbb{N}$. 
\end{definition}

\begin{definition} \quad  Let $N = \{N(u_{i})\}^k_{i=1}$ be a nbh sequence  in a graph $G$. Then the \emph{nbh sequence graph of} $N$ in $G$ denoted by $NSG[N]_G$ or simply $NSG[N]$ and is defined as the union of all induced subgraphs of closed nbh $N[u_{i}]$ in $G$, $i$ = 1 to $k$. i.e., $NSG[N]$ = $\bigcup^k_{i=1} {<N[u_{i}]>}$ which is a subgraph of $G$, $k\in\mathbb{N}$. 
\end{definition}

\begin{exm} \label{ex1} \quad Consider grid graph $G = P_{5} \Box P_{4}$ whose vertices are labeled as given in Figure 1. See Figure 1. In $G$, 
	
	$N[u^{(2)}_1]$ = $\{u^{(2)}_1, u^{(1)}_1, u^{(3)}_1, u^{(2)}_2\}$, 
	
	$N[u^{(3)}_2]$ = $\{u^{(3)}_2, u^{(3)}_1, u^{(2)}_2, u^{(4)}_2, u^{(3)}_3\}$, 
	
	$N[u^{(4)}_3]$ = $\{u^{(4)}_3, u^{(4)}_2, u^{(3)}_3, u^{(4)}_4\}$,
	
	$N[u^{(2)}_4]$ = $\{u^{(2)}_4, u^{(2)}_3, u^{(1)}_4, u^{(3)}_4\}, u^{(2)}_5$,
	
	$N[u^{(3)}_5]$ = $\{u^{(3)}_5, u^{(3)}_4, u^{(2)}_5, u^{(4)}_5\}$.
	
	Let $NSG[N]_1$ = $<N[u^{(2)}_1]>$ $\cup$ $<N[u^{(3)}_2]>$  $\cup$ $<N[u^{(4)}_3]>$, 
	
	$NSG[N]_2$ = $<N[u^{(2)}_4]>$  $\cup$ $<N[u^{(3)}_5]>$ and 
	
	$NSG[N]_3$ = $<N[u^{(2)}_1]>$ $\cup$ $<N[u^{(3)}_2]>$ 
	
	\hfill  $\cup$ $<N[u^{(4)}_3]>$ $\cup$ $<N[u^{(2)}_4]>$  $\cup$ $<N[u^{(3)}_5]>$  = $NSG[N]_1$ $\cup$ $NSG[N]_2$. 
	
	Here, $NSG[N]_1$, $NSG[N]_2$ and $NSG[N]_3$ are nbh sequence graphs and these induced subgraphs can be identified with respective vertices with yellow and blue colors and edges with green color in $G$ and are given in Figure 1. 	Clearly, $NSG[N]_1$ and $NSG[N]_2$ are connected subgraphs of $G$ whereas $NSG[N]_3$ is a disconnected subgraph of $G$.
	
\vspace{0.2cm} 
\begin{center}
	\begin{tikzpicture}  
	[scale=.7,auto=center,every node/.style={draw,circle,fill=white!20}] 
		
	\node (1) at (-8.5,2)[circle,fill=blue!20] {$\tiny u_{1}^{(1)}$};
	\node (10) at (-8.5,-0.5)[circle,fill=yellow!50] {$ u_{1}^{(2)}$};
	\node (13) at (-8.5,-3) [circle,fill=blue!20]  {$u_{1}^{(3)}$};
	\node (16) at (-8.5,-5.5) {$ u_{1}^{(4)}$};
	
	\node (4) at (-6,2)  {$u_{2}^{(1)}$};
	\node (2) at (-6,-0.5)[circle,fill=blue!20] {$ u_{2}^{(2)}$};
	\node (11) at (-6,-3)[circle,fill=yellow!50]  {$ u_{2}^{(3)}$};
	\node (14) at (-6,-5.5) [circle,fill=blue!20]  {$ u_{2}^{(4)}$};
	
	\node (18) at (-3.5,2)  {$u_{3}^{(1)}$};	
	\node (5) at (-3.5,-0.5) [circle,fill=blue!20]  {$ u_{3}^{(2)}$};
	\node (7) at (-3.5,-3) [circle,fill=blue!20] {$u_{3}^{(3)}$};
	\node (12) at (-3.5,-5.5) [circle,fill=yellow!50]{$ u_{3}^{(4)}$};
	
	\node (19) at (-0.5,2) [circle,fill=blue!20]  {$u_{4}^{(1)}$};	
	\node (20) at (-0.5,-0.5) [circle,fill=yellow!50] {$ u_{4}^{(2)}$};
	\node (21) at (-0.5,-3) [circle,fill=blue!20]  {$u_{4}^{(3)}$};
	\node (22) at (-0.5,-5.5) [circle,fill=blue!20] {$ u_{4}^{(4)}$};
	
	\node (23) at (2,2)  {$u_{5}^{(1)}$};	
	\node (24) at (2,-0.5) [circle,fill=blue!20]  {$ u_{5}^{(2)}$};
	\node (25) at (2,-3) [circle,fill=yellow!50] {$u_{5}^{(3)}$};
	\node (26) at (2,-5.5) [circle,fill=blue!20]  {$ u_{5}^{(4)}$};
	
	\draw (1) -- (4);
	\draw  (1)[line width = 0.1cm,draw=green] -- (10);
	\draw (2) -- (4);
	\draw (2) -- (5);
	\draw  (2)[line width = 0.1cm,draw=green]  -- (10);
	\draw (2)[line width = 0.1cm,draw=green]  -- (11);
	\draw (7) -- (5);
	\draw (7)[line width = 0.1cm,draw=green] -- (11);
	\draw (7) -- (21);
	\draw (7)[line width = 0.1cm,draw=green]  -- (12);
	\draw (10)[line width = 0.1cm,draw=green]  -- (13);
	\draw (11)[line width = 0.1cm,draw=green] -- (13);
	\draw (11)[line width = 0.1cm,draw=green]  -- (14);
	
	\draw (12)[line width = 0.1cm,draw=green]  -- (22);
	\draw (13) -- (16);
	
	\draw (14) -- (16);
	\draw (12)[line width = 0.1cm,draw=green]  -- (14);
	\draw (4) -- (18);
	
	\draw (20)[line width = 0.1cm,draw=green]  -- (19);
	\draw (20)[line width = 0.1cm,draw=green]  -- (5);
	\draw (20)[line width = 0.1cm,draw=green]  -- (21);
	\draw (20)[line width = 0.1cm,draw=green]  -- (24);
	
	\draw (25)[line width = 0.1cm,draw=green]  -- (24);
	\draw (25)[line width = 0.1cm,draw=green]  -- (21);
	\draw (25)[line width = 0.1cm,draw=green]  -- (26);
	
	\draw (18) -- (5);				
	\draw (18) -- (19);
	\draw (21) -- (22);		
	\draw (19) -- (23);		
	\draw (22) -- (26);
	\draw (23) -- (24);

	\end{tikzpicture} \\ 
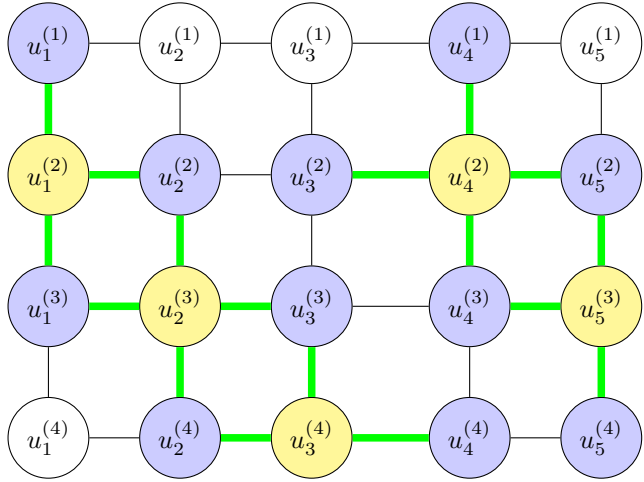
	\vspace{0.2cm} Figure 1. Grid graph $G = P_{5} \Box P_{4}$ \\
	\label{fig 7}
	
\end{center}		
\end{exm}

\begin{definition} \quad 	A nbh walk $N_1 N_2 \cdots  N_k$ is called a \emph{nbh trail} if every pair of non-consecutive terms(nbhs) of the nbh sequence has at the most one common vertex, except the first and the last terms. That is $|N_i \cap N_j|$ $\leq$ 1 for every $i,j$ $\ni$ $1 < |i-j| <$ $k-1$, $1 \leq i,j \leq k$. 
		
Thus, in a nbh trail $N_1 N_2 \cdots  N_k$, there is no restriction on the elements of $N_1$ and $N_k$ whereas $|N_i \cap N_j|$ $\leq$ 1 for every $i,j$ $\ni$ $1 < |i-j| < k-1$, $1 \leq i,j \leq k$. 
	\end{definition}

	\begin{definition} \quad	A nbh trial $N_1 N_2 \cdots  N_k$ is called a \emph{nbh path} or a \emph{nbh chain} if non-consecutive terms(nbhs) of the nbh sequence have no common vertex, except the first and the last terms. 
		
	Nbh chain $N_1 N_2 \cdots  N_k$ is said to be of \emph{length} $k$, $k\in\mathbb{N}$.				
\end{definition}				

\begin{definition} \quad	A nbh chain $N_1 N_2 \cdots N_k$ is called a \emph{closed nbh chain} or \emph{closed nbh path} or \emph{nbh cycle} if $N_1 \cap N_k \neq \emptyset$ and $N_1 \neq N_k$ when $k = 2$ and $N_1 \cap N_k \neq \emptyset$ when $k \geq 3$. Otherwise, it is called an \emph{open nbh chain}.
\end{definition}

\begin{definition} \quad  A nbh sequence $\{N(u^{(j)}_{i})\}^k_{j=1}$ in a graph $G$ is called a  \emph{connected nbh sequence} if between any two distinct nbhs(terms) of the nbh sequence, there exists at least one nbh path, $u^{(j)}_{i}\in V(G)$, $k \geq 2$, $1 \leq j \leq k$, $i\in\mathbb{N}$. Otherwise, it is called a \emph{disconnected nbh sequence} in $G$. 
	
{\rm Thus, in graph $G$ as given in Example \ref{ex1}, $NSG[N]_1$ and $NSG[N]_2$ are connected nbh sequences whereas $NSG[N]_3$ is a disconnected nbh sequence.}	
\end{definition}
	
\begin{exm} \label{ex2} \quad Let $G$ be a graph and $N_{i}$ = $N(u_i)$, $u_i\in V(G)$ and $1 \leq i \leq 11$. In Figure 2, each colored region represents a nbh $N_i$, $1 \leq i \leq 11$; any two regions are overlapping if and only if their corresponding nbhs have common element(s)(vertex or vertices). In each of the figures, Figures 2, 2(a), 2(b) and 3, each region shown corresponds to a given nbh in the graph $G$ and not its corresponding subgraph or induced subgraph in $G$. To understand more about these definitions on nbh sequences, we present these examples.
		
	
	\begin{center}
		\begin{tikzpicture}
			\draw[fill=green!70,fill opacity=0.6,rotate=0](2,-1.1) ellipse (0.4cm and 1.2cm);
			\draw[fill=magenta!70,fill opacity=0.8,rotate=-35] (1.7,-0.6) ellipse (1.5cm and 0.5cm);
			\draw[fill=red!70,fill opacity=0.6,rotate=0] (0.15,0) ellipse (0.6cm and 1.7cm);
			\draw[fill=brown!70,fill opacity=0.4,rotate=30] (1.8,0.6) ellipse (1.5cm and 0.5cm);
			\draw[fill=orange!70,fill opacity=0.4] (3.5,2.1) ellipse (1.5cm and 0.6cm);
			\draw[fill=red!70,fill opacity=0.5,rotate=35] (5,-2) ellipse (0.6cm and 1.5cm);
			\draw[fill=green!70,fill opacity=0.8] (4.7,0) ellipse (1.5cm and .8cm);
			\draw[fill=blue!70,fill opacity=0.8] (2.5,0) ellipse (2cm and .5cm);
			\draw[fill=yellow!70,fill opacity=0.5] (0,0) ellipse (1.5cm and .8cm);
			\draw[fill=pink!70,fill opacity=0.5] (-2.5,0) ellipse (1.5cm and .8cm);
			\draw[fill=magenta!70,fill opacity=0.8,rotate=-35] (5.7,2.5) ellipse (1.5cm and .5cm);
			\node at (-3,0) {$N_{1}$};
			\node at (-0.7,0) {$N_{2}$};
			\node at (2.6,0) {$N_{3}$};
			\node at (5,0) {$ N_{4}$};
			\node at (5,1.3) {$ N_{5}$};
			\node at (3.5,2) {$ N_{6}$};
			\node at (1,1.5) {$ N_{7}$};
			\node at (0.1,0) {$ N_{8}$};
			\node at (1.2,-1.5) {$ N_{9}$};
			\node at (6,-1.3) {$ N_{10}$};
			\node at (2,-1) {$ N_{11}$};
		\end{tikzpicture}
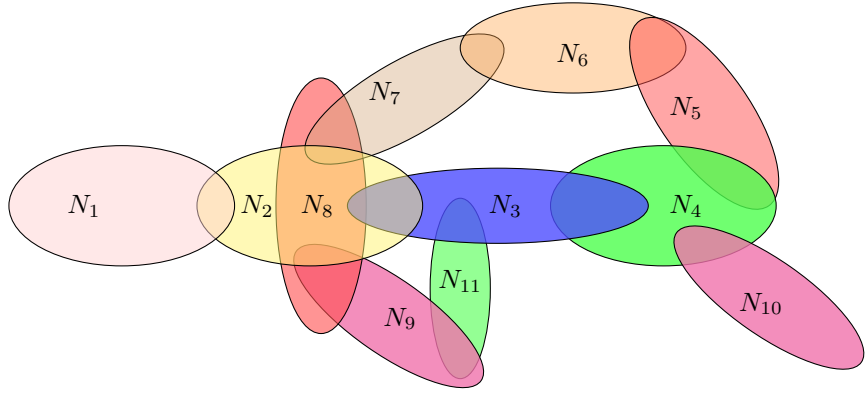
	 $Figure ~2.$ Nbh sequence graph of $N_1 N_2 . . . N_{11}$ \label{fig 2}
	\end{center}
	
	Consider the following sequence of neighbourhoods of vertices in $G$. Let them be
	
	NS1 = $N_1 N_2 N_3 N_4 N_5 N_6 N_7 N_8 N_9 $, 
	
	NS2 = $N_1 N_2 N_8 N_7 N_6 N_5 N_4 N_{10}$, 
	
	NS3 = $N_1 N_2 N_3 N_4 N_5 N_6 N_7 N_8 N_3 N_4 N_{10}$, 
	
	NS4 = $N_1 N_2 N_3 N_4 N_5 N_6 N_7 N_8 N_2 N_3 N_4 N_{10}$, 
	
	NS5 = $N_1 N_2 N_3 N_4 N_5 N_6 N_5 N_4 N_{10}$,
	
	NS6 = $N_1 N_2 N_3 N_2 N_1 N_2 N_9$, 
	
	NS7 = $N_9 N_8 N_2 N_1 N_2 N_9$, 
	
	NS8 = $N_9 N_8 N_3 N_{11} N_9$ = $N_9 N_8 N_3 N_{11}$, 
	
	NS9 = $N_8 N_3 N_4 N_5 N_6 N_7 N_8$ = $N_8 N_3 N_4 N_5 N_6 N_7$,
	
	NS10 = $N_2 N_3 N_4 N_5 N_6 N_7 N_2$ = $N_2 N_3 N_4 N_5 N_6 N_7$,
	
	NS11 = $ N_2 N_3 N_4 N_5 N_6 N_7 N_8$, 
	
	NS12 = $N_1 N_2 N_3 N_4 N_5 N_6 N_7 N_8 N_9 N_{11} N_3 N_4 N_3 N_8 N_2 N_1$, see Figure 2(a), 
	
	NS13 = $N_8 N_9 N_{11} N_3 N_4 N_5 N_6 N_7$, see Figure 2(b), 
	
	NS14 = $N_8 N_9 N_{11} N_3$
	
	NS15 = $N_2 N_3 N_4 N_5 N_6 N_7 N_9$,
	
	NS16 = $N_1 N_8 N_3$, 
	
	NS17 = $N_9 N_8 N_6 N_5$,
	
	NS18 = $N_8 N_9 N_{10} N_{11} N_3 N_4 N_{10}$. 
	\\	
	Now, let us check the above sequences of nbhs for nbh walk, nbh trail, nbh path/chain, nbh cycle.
	
	\begin{enumerate}	
		\item [\rm (i)] {\bf Nbh Walks:}\quad 
		
		NS1 to NS14 are nbh walks.
		
		NS15 to NS18 are not nbh walks since 
		
		(a) in NS15, consecutive terms(nbhs) $N_7$ and $N_9$ are disjoint;
		
		(b) in NS16, consecutive terms(nbhs) $N_1$ and $N_8$ are disjoint;
		
		(c) in NS17,  consecutive terms(nbhs) $N_8$ and $N_6$ are disjoint;
		
		(d) in NS18,  consecutive pair of terms(nbhs) $N_9$ and $N_{10}$ 
		
		\hfill (as well as $N_{10}$ and $N_{11}$) are disjoint. 
		
	\item [\rm (ii)] {\bf Nbh Trails:}\quad  Among the nbh walks,
	
	NS8, NS9, NS10, NS13, NS14 are nbh chains/paths;
	
	NS1 is a nbh trail 
	
	\hfill if $|N_7 \cap N_2| = 1$, $|N_8 \cap N_2| = 1$, $|N_8 \cap N_3|$ = 1 and $|N_9 \cap N_2| = 1$; 
	
	NS1 is not a nbh trail 
	
	\hfill if $|N_7 \cap N_2| > 1$ or $|N_8 \cap N_2| > 1$ or $|N_8 \cap N_3| > 1$ or $|N_9 \cap N_2| > 1$;
	
	NS2 is a nbh trail  if $|N_7 \cap N_2| = 1$ and is not a nbh trail  if $|N_7 \cap N_2| > 1$;
	
	NS3 is not a nbh trail since $N_3$ (and $N_4$) occurs twice in the nbh sequence;
	
	NS4 is not a nbh trail since $N_3$ (as well as $N_4$) occurs twice in NS4;
	
	NS5 is not a nbh trail since $N_4$ (as well as $N_5$) occurs twice in NS5;

	NS6 is not a nbh trail 
	
	\hfill since $N_1$ (as well as $N_2$) occurs more than once in NS6;
	
	NS7 is not a nbh trail since $N_2$ occurs more than once in NS7;

	NS11 is a nbh trail if $|N_7 \cap N_2| = 1$ and $|N_8 \cap N_3| = 1$;

	NS11 is not a nbh trail if $|N_7 \cap N_2| > 1$ or $|N_8 \cap N_3| > 1$;

	NS12 is not a nbh trail since $N_2$ (and $N_3$, $N_8$) occurs more than once in
	
	\hfill  NS12 and is (are) not an end term/nbh of NS12 (See Figure 2(a).).

	\item [\rm (iii)] {\bf Nbh Paths:}\quad  
	\\
	Among the nbh trails,
	
	NS2, NS8, NS9, NS10, NS13, NS14 are nbh chains/paths;
	
	NS1 is not a nbh path since 
	
	\hfill $|N_7 \cap N_2| \neq \emptyset$ (and $|N_8 \cap N_2|  \neq \emptyset$, $|N_8 \cap N_3|  \neq \emptyset$, $|N_9 \cap N_2|  \neq \emptyset$);
	
	NS3 is not a nbh path since $|N_8 \cap N_2| \neq \emptyset$;
	
	NS4 is not a nbh path since $N_3$ (as well as $N_4$) 
	
	\hfill is not an end term/nbh in NS4 and occurs twice in NS4;
	
	NS5 is not a nbh path since $N_4$ (as well as $N_5$) 
	
	\hfill is not an end term/nbh in NS5 and occurs twice in NS5;
	
	NS6 is not a nbh path since $N_1$ (as well as $N_2$) occurs
	
	\hfill  more than once in NS6 and is not an end term/nbh of NS7;
	
	NS7 is not a nbh path since $N_2$ occurs 
	
		\hfill more than once and is not an end term/nbh of NS7;
	
	NS11 is not a nbh path since $N_2$ and $N_7$ 
	
		\hfill are non-adjacent terms/nbhs but  $N_{2} \cap N_{7} \neq \emptyset$;

    NS12 is not a nbh path since $N_3$ occurs more than once in NS12
    
    \hfill   and is not an end term/nbh of NS12 (See Figure 2(a).);
    
    NS13 is a nbh trail and also a nbh path since
    
    \hfill   $N_8 N_9 N_{11} N_3 N_4 N_5 N_6$ is a nbh path  and $N_7$ is adjacent only to 
    
    \hfill  nonadjacent terms $N_6$ and $N_8$ in NS13 (See Figure 2(b).).

	\item [\rm (iv)] {\bf Nbh Cycles:}\quad  
	
	NS1 is not a nbh cycle since $N_{1} \cap N_{9} = \emptyset$  but it contains the nbh cycle $N_3 N_4 N_5 N_6 N_7 N_8$. 
	
	Among the nbh paths NS2, NS8, NS9, NS10, NS13 and NS14, 
	
	NS8, NS9, NS10 and NS14 are nbh cycles whereas 
	
	NS2 and NS13 are not nbh cycles since
	
	(a) in NS2, $N_{1} \cap N_{10} = \emptyset$;
	
	(b) NS13 contains nbh cycles $N_8 N_9 N_{11} N_3$ and $N_8 N_3 N_4 N_5 N_6 N_7$ 
	
	 \hfill  (more than one cycle in NS13).	
		
	\end{enumerate} 	
\begin{center}
	\begin{tikzpicture}
	\draw[fill=magenta!70,fill opacity=0.8,rotate=-35] (1.7,-0.4) ellipse (1.3cm and 0.4cm);
	\draw[fill=red!70,fill opacity=0.6,rotate=0] (0.15,0) ellipse (0.6cm and 1cm);
	\draw[fill=brown!70,fill opacity=0.4,rotate=35] (1.6,0.4) ellipse (1.2cm and 0.5cm);
	\draw[fill=orange!70,fill opacity=0.4] (2.6,1.5) ellipse (1cm and 0.5cm);
	\draw[fill=red!70,fill opacity=0.5,rotate=35] (3.6,-1.8) ellipse (0.5cm and 1cm);
	\draw[fill=green!70,fill opacity=0.8] (3.5,0) ellipse (1cm and .5cm);
	\draw[fill=blue!70,fill opacity=0.8] (1.6,0) ellipse (1.2cm and .4cm);
	\draw[fill=yellow!70,fill opacity=0.5] (0,0) ellipse (1.1cm and .6cm);
	\draw[fill=pink!70,fill opacity=0.5] (-1.9,0) ellipse (1cm and .4cm);
	\draw[fill=green!70,fill opacity=0.6,rotate=0](2,-1.1) ellipse (0.4cm and 1cm);
	
	\node at (-2,0) {$N_{1}$};
	\node at (-0.7,0) {$N_{2}$};
	\node at (1.6,0) {$N_{3}$};
	\node at (3.5,0) {$ N_{4}$};
	\node at (3.8,1) {$ N_{5}$};
	\node at (2.8,1.5) {$ N_{6}$};
	\node at (1,1.3) {$ N_{7}$};
	\node at (0.1,0) {$ N_{8}$};
	\node at (1.2,-1.5) {$ N_{9}$};
	\node at (2,-1) {$ N_{11}$};
	\end{tikzpicture}
	
	\vspace{0.2cm} $Figure ~2(a).$ $NS12 $ \label{fig 2(a)}\\
	
	\vspace{0.5cm}
	\begin{tikzpicture}
	\draw[fill=magenta!70,fill opacity=0.8,rotate=-35] (1.7,-0.4) ellipse (1.3cm and 0.4cm);
	\draw[fill=red!70,fill opacity=0.6,rotate=0] (0.15,0) ellipse (0.6cm and 1cm);
	\draw[fill=brown!70,fill opacity=0.4,rotate=35] (1.6,0.4) ellipse (1.2cm and 0.5cm);
	\draw[fill=orange!70,fill opacity=0.4] (2.6,1.5) ellipse (1cm and 0.5cm);
	\draw[fill=red!70,fill opacity=0.5,rotate=35] (3.6,-1.8) ellipse (0.5cm and 1cm);
	\draw[fill=green!70,fill opacity=0.8] (3.5,0) ellipse (1cm and .5cm);
	\draw[fill=blue!70,fill opacity=0.8] (1.6,0) ellipse (1.2cm and .4cm);
	\draw[fill=green!70,fill opacity=0.6,rotate=0](2,-1.1) ellipse (0.4cm and 1cm);
	
	\node at (1.6,0) {$N_{3}$};
	\node at (3.5,0) {$ N_{4}$};
	\node at (3.8,1) {$ N_{5}$};
	\node at (2.8,1.5) {$ N_{6}$};
	\node at (1,1.3) {$ N_{7}$};
	\node at (0.1,0) {$ N_{8}$};
	\node at (1.2,-1.5) {$ N_{9}$};
	\node at (2,-1) {$ N_{11}$};
	\end{tikzpicture}\\
	\vspace{0.2cm} $Figure~ 2(b).$ $ NS13 $ \label{fig 2(b)}
\end{center}
\end{exm}
	\begin{exm} \quad Let $G$ be a graph and $N_{i}$ = $N(u_i)$, $u_i\in V(G)$ and $12 \leq i \leq 22$. In Figure 3, each colored region represents a nbh $N_i$, $12 \leq i \leq 22$; any two regions are overlapping if and only if their corresponding nbhs have common element(s)(vertex or vertices).	Consider the following nbh sequences. Let 
		
	$NS19$ =	$N_{12} N_{13} N_{14} N_{15} N_{16} N_{17} N_{18} N_{19} N_{20} N_{21}$; 
		
	$NS20$ = $N_{13} N_{14} N_{15} N_{16} N_{17} N_{18} N_{19}$;
		
	$NS21$ = $N_{12} N_{13} N_{14} N_{15} N_{16} N_{17} N_{18} N_{19}$ and 
		
	$NS22$ = $N_{13} N_{14} N_{15} N_{16} N_{17} N_{18} N_{19} N_{13} N_{20} N_{21}$.
	
	By checking these nbh sequences, we identify them as follows.
	
	(a) $NS19$ is not a nbh walk since $N_{19} \cap N_{20} = \emptyset$ eventhough
	
	\hfill  $N_{19}$ and $N_{20}$ are consecutive terms/nbhs in NS19. 
	
	(b) NS20  is a nbh cycle since NS20  is a nbh chain and $N_{13} \cap N_{19} \neq \emptyset$. 
	
	(c) $NS21$ is a nbh chain but not nbh cycles since $N_{12} \cap N_{19} = \emptyset$.
	
	(d) $NS22$ is a nbh chain but not nbh cycles since $N_{13} \cap N_{21} = \emptyset$.
	\end{exm}
	\begin{center}
		\begin{tikzpicture}
			\draw[fill=black!20,fill opacity=0.8] (2.8,-1.7) ellipse (1.5cm and 0.5cm);	
			\draw[fill=magenta!70,fill opacity=0.8,rotate=-50] (1.6,-0.05) ellipse (1.2cm and 0.5cm);
			\draw[fill=red!40,fill opacity=0.6,rotate=-20] (0,1.3) ellipse (0.5cm and 1.1cm);
			\draw[fill=brown!70,fill opacity=0.4] (1.8,2) ellipse (1.4cm and 0.5cm);
			\draw[fill=orange!70,fill opacity=0.4] (3.7,2) ellipse (1.5cm and 0.6cm);
			\draw[fill=red!70,fill opacity=0.5,rotate=30] (5,-1.5) ellipse (0.6cm and 1.5cm);
			\draw[fill=green!70,fill opacity=0.8] (4.5,0) ellipse (1.5cm and .6cm);
			\draw[fill=blue!70,fill opacity=0.8] (2.5,0) ellipse (1.5cm and .6cm);
			\draw[fill=yellow!70,fill opacity=0.5] (0,0) ellipse (1.5cm and .6cm);
			\draw[fill=pink!70,fill opacity=0.5] (-2.5,0) ellipse (1.5cm and .6cm);
				\draw[fill=red!40,fill opacity=0.6,rotate=-30] (4.5,1.2) ellipse (0.5cm and 1.3cm);
				\draw[fill=red!40,fill opacity=0.6,rotate=-30] (2.2,2.5) ellipse (0.4cm and 1cm);
			\node at (-3,0) {$N_{12}$};
			\node at (-0.7,0) {$N_{13}$};
			\node at (2.5,0) {$N_{14}$};
			\node at (4.3,0) {$ N_{15}$};
			\node at (5,1.3) {$ N_{16}$};
			\node at (3.8,2.1) {$ N_{17}$};
			\node at (1.8,2) {$ N_{18}$};
			\node at (0.2,1) {$ N_{19}$};
			\node at (0.9,-1.2) {$ N_{20}$};
			\node at (2.5,-1.8) {$N_{21}$};
			\node at (4.5,-1) {$ N_{22}$};
			\node at (3.1,1) {$ N_{23}$};
		\end{tikzpicture}
		
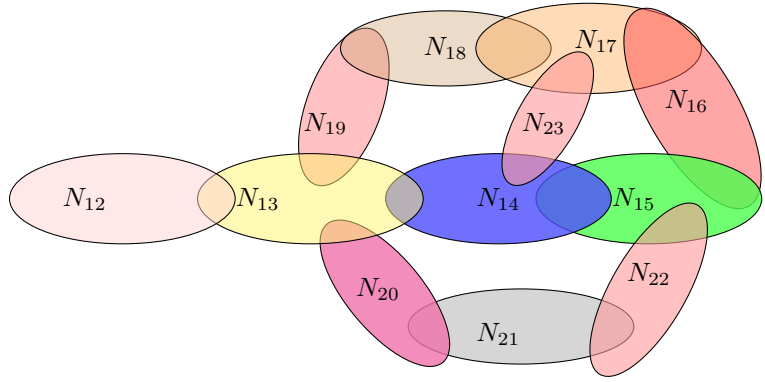
		\vspace{0.2cm} $Figure ~3.$ Nbh sequence graph of $N_{12} N_{13} . . . N_{22}$ \label{fig 3}
	\end{center}
		
	\section{Nbh chains of Type-1 (NC-T1) and Type-2 (NC-T2) on NDM graphs}
	
	 In this section, we define nbh chains of Type-1 (NC-T1) and Type-2 (NC-T2) and prove that (i) for $k \geq 2$ and $n \geq 3$, graph $P_{k} \Box C_{n}$ contains NC-T2, $n,k \in \mathbb{N}$; (ii) graph containing NC-T1  of even length is NDM and (iii) for $m \geq 2$ and $n \geq 3$,  cylindrical grid graphs $P_{2m} \Box C_n$ are NDM,  $m,n\in\mathbb{N}$.
	 	
	\begin{definition} \quad {\rm[\textit{Neighbourhood Chain of Type 1 (NC-T1)}]}
		
		Let $n \geq 2 $,  $m \geq n+1$ and $G$ be a graph, $m,n \in \mathbb{N}$. Let $\{u_{1},u_{2},\ldots,u_{n}\}$ and $\{v_{1},v_{2},\ldots,v_{m}\}$ be two disjoint subsets of $V(G)$, $u_{1} \neq u_{n}$ and $N_{i}$ = $N(u_i)$ = $\{v_j\in V(G): u_{i}v_j\in E(G)$, $ 1\leq j\leq m \}$ $\neq$ $\emptyset$ for $i = 1,2,\ldots,n$. Let $NSG$ = $N_1 N_2 \cdots N_{n}$ be a neighbourhood chain in $G$. Then the nbh chain $NSG$ is said to be of Type 1 (NC-T1) if it satisfies the following properties:
	\begin{enumerate}
		\item [\rm (i) ] $N_{1}\setminus N_{2} = \left\{v_{1}\right\}$ and $N_{n} \setminus N_{n-1}$ = $\left\{v_{n}\right\}$;
		
		\item [\rm (ii) ]  $N_{i+1} \setminus N_{i} \subset N_{i+2}$, $i = 1,2,\ldots,n-2$ and $\vert N_{i} \setminus N_{i+1} \vert, \vert N_{i+1} \setminus N_{i} \vert$ $\geq 1$ for $i = 1,2,\ldots,n-1$. 
	\end{enumerate}	
	\end{definition}

	\begin{note}\quad {\rm The conditions given in the above definition of NC-T1 are a more generalization of conditions given in Theorem \ref{a2}.}
\end{note}

\begin{definition} \quad \label{d3} {\rm[\textit{Type 2 Neighbourhood Chains (NC2-T2)}]}
	
	Let $G$ be a graph and $N_1 N_2 \cdots N_{n}$ and  $N_{n+1} N_{n+2} \cdots N_{2n}$ be two nbh chains of Type-1 (NC-T1) in $G$, $n \geq 2$ and $ n\in \mathbb{N}$. Then the two nbh chains of Type-1 are called as a Type 2 Neighbourhood Chains (NC-T2) if they satisfy the following properties.
	\begin{enumerate}
		\item [\rm (i) ]  $N_{1} \cap N_{2} \cap N_{n+1} = N_{n+1} \setminus N_{n+2} \neq \emptyset$ and $ N_{n} \cap N_{2n} \cap N_{2n-1} = N_{n} \setminus N_{n-1} \neq \emptyset $;
	\item [\rm (ii) ]  $ N_i \cap N_{i+1} \cap N_{n+i} \cap N_{n+i-1} \neq \emptyset$ for $i = 2,3,\ldots,n-1$.
	\end{enumerate}
\end{definition}

Our next theorem is to prove that graph $P_{k} \Box C_{n}$ contains a Type-2 nbh chains for $k \geq 2$, $n \geq 3$ and $n,k \in \mathbb{N}$. Before proving this theorem, we present two examples to illustrate this result.

\begin{exm} \quad In $P_{4} \Box C_{3}$, the nbh chain  $N_{1}^{(1)} N_{2}^{(2)} N_{3}^{(3)} N_{4}^{(1)}$, as given in Figure 4, is of Type 1 where $N_i^{(j)}$ = $N(u_i^{(j)})$, $1 \leq i \leq 4$ and $1 \leq j \leq 3$. See Figure 4.
\end{exm}	
	In $P_{4} \Box C_{3}$, vertices are labeled as given in Figure 4. Let $u_1 = u_{1}^{(1)}$, $u_2 = u_{2}^{(2)}$, $u_3 = u_{3}^{(3)}$, $u_4 = u_{4}^{(1)}$, $v_1 = u_{1}^{(3)}$, $v_2 = u_{2}^{(1)}$, $v_3 = u_{3}^{(2)}$, $v_4 = u_{4}^{(2)}$, $v_5 = u_{4}^{(3)}$, $v_6 = u_{1}^{(2)}$, $v_7 = u_{2}^{(3)}$, $v_8 = u_{3}^{(1)}$. Then,
	
\begin{align*}
N_{1}^{(1)}~&=~N(u_{1}^{(1)})~=~N(u_1)~=~\{v_1,v_2,v_6\}, \\
N_{2}^{(2)}~&=~N( u_{2}^{(2)})~=~N(u_2)~=~\{v_2,v_6,v_3,v_7\}, \\
N_{3}^{(3)}~&=~N(u_{3}^{(3)})~=~N(u_3)~=~\{v_3,v_7,v_5,v_8\},\\
N_{4}^{(1)}~&=~N(u_{4}^{(1)})~=~N(u_4)~=~\{v_4,v_5,v_8\}.
\end{align*}
\begin{align*}
N_{1}^{(1)} \cap  N_{2}^{(2)} = \{v_2,v_6\} \neq \emptyset,\\
N_{2}^{(2)} \cap N_{3}^{(3)} = \{v_3,v_7\} \neq \emptyset,\\
N_{3}^{(3)} \cap N_{4}^{(1)} = \{v_5,v_8\} \neq \emptyset.
\end{align*}		
	\begin{center}
		\begin{tikzpicture}
			[node distance={20mm}, thick, main/.style = {draw,circle,fill=blue!20}]]

			\node[main] (2) [circle,fill=red!50][label=180:$u_1$]{\tiny $=u_{1}^{(1)}$}; 
			\node[main] (1) [below left of=2] [label=180:$v_1$]{\tiny $=u_{1}^{(3)}$}; 
			\node[main] (3) [below right of=1][label=180:$v_6$] {\tiny $=u_{1}^{(2)}$}; 
			\path (2) edge[left=60] (1);
			\path (2) edge[left=60] (3);
			\path (1) edge[ left=60] (3);

			\node[main] (8) [above right of=3][label=175:$v_7$]{\tiny $=u_{2}^{(3)}$};
			\node[main] (9) [above right of=8] [label=175:$v_2$]{\tiny $=u_{2}^{(1)}$};
			\node[main] (10) [below right of=8,][circle,fill=red!50][label=190:$u_2$] {\tiny $ =u_{2}^{(2)}$};
			\path (8) edge[left=30] (9);
			\path (8) edge[left=30] (10);
			\path (9) edge[left=30] (10);
			\path (9) edge[left=30] (2);
			\path (1) edge[left=30] (8);
			
			\path (3) edge[left=30] (10);
			
			\node[main] (15) [above right of=10][circle,fill=red!50][label=175:$u_3$]{\tiny $ =u_{3}^{(3)}$};
			\node[main] (16) [above right of=15][label=175:$v_8$] {\tiny $=u_{3}^{(1)}$};
			\node[main] (17) [below right of=15,][label=190:$v_3$] {\tiny $ =u_{3}^{(2)}$};
			\path (8) edge[left=30] (15);	
			\path (15) edge[left=30] (16);
			\path (15) edge[left=30] (17);
			\path (16) edge[left=30] (17);
			\path (9) edge[left=30] (16);
			\path (10) edge[left=30] (17);
			
			\node[main] (22) [above right of=17][label=0:$v_5$]{\tiny $ u_{4}^{(3)}=$};
			\node[main] (23) [above right of=22][circle,fill=red!50][label=0:$u_4$] {\tiny $ u_{4}^{(1)}=$};
			\node[main] (24) [below right of=22,][label=0:$v_4$] {\tiny $ u_{4}^{(2)}=$};
			\path (16) edge[left=30] (23);
			\path (22) edge[left=30] (23);
			\path (22) edge[left=30] (24);
			\path (23) edge[left=30] (24);
			\path (15) edge[left=30] (22);
			\path (17) edge[left=30] (24);	\end{tikzpicture}\\
		\vspace{0.2cm} Figure 4.~ $P_4 \Box C_3$
	\end{center}

\vspace{0.2cm}
	Thus, $N_{1}^{(1)} N_{2}^{(2)} N_{3}^{(3)} N_{4}^{(1)}$ is a nbh chain in $P_{4} \Box C_{3}$. Also,
	\begin{align}
		N_{1}^{(1)} \setminus  N_{2}^{(2)} ~&=~\{v_1\} ~\& ~N_{4}^{(1)} \setminus N_{3}^{(3)} ~=~\{v_4\}\label{b};\\
		 N_{2}^{(2)} \setminus N_{1}^{(1)} ~&=~\{v_3,v_7\}~\subset ~ \{v_3,v_7,v_5,v_8\}~= N_{3}^{(3)};\label{c}\\
		N_{3}^{(3)} \setminus  N_{2}^{(2)} ~&=~\{v_5,v_8\}~\subset \{v_4,v_5,v_8\} = N_{4}^{(1)};\label{d}\\
		\vert  N_{2}^{(2)} \setminus N_{1}^{(1)} \vert &~=~2 >~1;~ ~\vert N_{3}^{(3)} \setminus  N_{2}^{(2)} \vert~=~2 >~1.\label{e} 		
	\end{align}
	
	This implies that $N_{1}^{(1)} N_{2}^{(2)} N_{3}^{(3)} N_{4}^{(1)}$ is a nbh chain of Type-1 in $P_4 \Box C_3$. \hfill $\Box$
	
	\vspace{.2cm}
	In Figure 5, we emphasis nbh chain $N_{1}^{(1)} N_{2}^{(2)} N_{3}^{(3)} N_{4}^{(1)}$ of Type-1 contained in $P_4 \Box C_3$. And for more clarity, we draw edges of $P_4 \Box C_3$ as continuous lines as well as dotted lines. Edges joining a point and its adjacent points in each nbh of the nbh chain $N_{1}^{(1)} N_{2}^{(2)} N_{3}^{(3)} N_{4}^{(1)}$ of Type-1 are drawn as continuous lines and all other edges of $P_4 \Box C_3$ are represented by dotted lines.\hfill $\Box$
	
		\begin{center} 
		\begin{tikzpicture}  
			[scale=.8,auto=center,every node/.style={draw,circle,fill=blue!20}] 
			
			\node (1) at (-9,2)[circle,fill=red!50][label=0:$u_{1}^{(1)}$]{$u_{1}=$};
			
			\node (2) at (-5.5,2)[circle,fill=red!50][label=0:$\tiny u_{2}^{(2)}$]{$u_{2}=$};
			\node (3) at (-10,4)  [label=0:$\tiny u_{1}^{(3)}$]{$v_{1}=$};
			\node (4) at (-7.3,4) [label=0:$\tiny u_{2}^{(1)}$]{$v_{2}=$};
			\node (5) at (-3.7,4) [label=0:$\tiny u_{3}^{(2)}$]{$v_{3}=$};
			\node (6) at (2,0) [label=0:$\tiny u_{4}^{(2)}$]{$v_{4}=$};
			\node (7) at (-2,2) [circle,fill=red!50][label=0:$\tiny u_{3}^{(3)}$]{$u_{3}=$};
			
			\node (8) at (1,2) [circle,fill=red!50][label=0:$\tiny u_{4}^{(1)}$]{$u_{4}=$};
			\node (9) at (-0.5,4) [label=0:$\tiny u_{4}^{(3)}$]{$v_{5}=$};
			\node (10) at (-7.3,0) [label=0:$\tiny u_{1}^{(2)}$]{$v_{6}=$};
			\node (11) at (-3.7,0) [label=0:$\tiny u_{2}^{(3)}$]{$v_{7}=$};
			\node (12) at (-0.5,0) [label=0:$\tiny u_{3}^{(1)}$]{$v_{8}=$};

			\draw (1) -- (3);
			\draw (1) -- (4);
			\draw (1) -- (10);
			\draw (2) -- (4);
			\draw (2) -- (5);
			\draw (2) -- (10);
			\draw (2) -- (11);
			\draw (7) -- (5);
			\draw (7) -- (11);
			\draw (7) -- (9);
			\draw (7) -- (12);
			\draw (8) -- (6);
			\draw (8) -- (9);
			\draw (8) -- (12);
			\draw[dashed]
			(3) to[out=-120,in=180] (10);
			\draw[dashed]
			(3) to[out=-15,in=180] (11);
			\draw[dashed]
			(4) to[out=-80,in=160] (11);
			\draw[dashed]
			(4) to[out=-10,in=160] (12);
			\draw[dashed]
			(5) to[out=-80,in=150] (12);
			\draw[dashed]
			(5) to[out=-15,in=135] (6);
			\draw[dashed]
			(9) to[out=-10,in=90] (6);
			
		\end{tikzpicture}\\  
		\vspace{0.2cm} Figure 5. $P_4 \Box C_3 $
	\end{center}

		\begin{exm}\quad Graph $P_{3} \Box C_{5}$ contains nbh chains of Type 2. See Figure 6. 
		\end{exm}
	In $P_{3} \Box C_{5}$, vertices are labeled as given in Figure 6. Let
			
	$u_1 = u_{1}^{(1)}$, $u_2 = u_{2}^{(2)}$, $u_3 = u_{3}^{(3)}$, $u_4 = u_{1}^{(3)}$, $u_5 = u_{2}^{(4)}$, $u_6 = u_{3}^{(5)}$, $v_1 = u_{1}^{(5)}$, $v_2 = u_{2}^{(1)}$, $v_3 = u_{3}^{(4)}$, $v_4 = u_{3}^{(2)}$, $v_5 = u_{1}^{(2)}$, $v_6 = u_{3}^{(1)}$, $v_7 = u_{2}^{(3)}$, $v_8 = u_{1}^{(4)}$, $v_9 = u_{2}^{(5)}$.
	
	At first, we show that $N_{1}^{(1)} N_{2}^{(2)} N_{3}^{(3)}$ and $N_{1}^{(3)} N_{2}^{(4)} N_{3}^{(5)}$ are nbh chains of Type 1 in $P_{3} \Box C_{5}$ where $N_i^{(j)}$ = $N(u_i^{(j)})$, $1 \leq i \leq 3$ and $1 \leq j \leq 5$. See Figure 6.
		
	\begin{center}
		\begin{tikzpicture}  
			[scale=.6,auto=center,every node/.style={draw,circle,fill=blue!20}] 
			
			\node (1) at (-10,3)[circle,fill=red!50][label=180:$u_{1}$]  {$\tiny =u_{1}^{(1)}$};
			\node (10) at (-8.5,-0.5)[label=180:$v_{5} $] {$ =u_{1}^{(2)}$};
			\node (13) at (-10,-4) [circle,fill=red!50] [label=180:$u_{4}$] {$=u_{1}^{(3)}$};
			\node (16) at (-11.5,-2) [label=180:$v_{8}$] {$ =u_{1}^{(4)}$};
			\node (3) at (-11.5,1)  [label=180:$v_{1}$]{$=u_{1}^{(5)}$};
						
			\node (4) at (-4,3)  [label=175:$ v_{2}$]{$=u_{2}^{(1)}$};
			\node (2) at (-2.5,-0.5)[circle,fill=red!50][label=10:$u_{2} $] {$ u_{2}^{(2)}=$};
			\node (11) at (-4,-4) [label=350:$v_{7}$] {$ u_{2}^{(3)}=$};
			\node (14) at (-5.5,-2) [circle,fill=red!50] [label=10:$u_{5}$] {$u_{2}^{(4)}=$};
			\node (17) at (-5.5,1) [label=175:$v_{9}$] {$ =u_{2}^{(5)}$};
				
			\node (18) at (2,3)  [label=0:$v_{6}$]{$u_{3}^{(1)}=$};	
			\node (5) at (3.5,-0.5) [label=0:$v_{4}$] {$ u_{3}^{(2)}=$};
			\node (7) at (2,-4) [circle,fill=red!50] [label=0:$ u_{3}$]{$u_{3}^{(3)}=$};
			\node (12) at (0.5,-2) [label=0:$ v_{3}$]{$ u_{3}^{(4)}=$};
			\node (15) at (0.5,1)  [circle,fill=red!50] [label=0:$u_{6}$]{$u_{3}^{(5)}=$};
			
			\draw (1) -- (3);
			\draw (1) -- (4);
			\draw (1) -- (10);
			\draw (2) -- (4);
			\draw (2) -- (5);
			\draw (2) -- (10);
			\draw (2) -- (11);
			\draw (7) -- (5);
			\draw (7) -- (11);
			
			\draw (7) -- (12);
			\draw (10) -- (13);
			\draw (11) -- (13);
			\draw (11) -- (14);
			\draw (12) -- (15);
			\draw (13) -- (16);
			\draw (14) -- (17);
			\draw (15) -- (18);	
			\draw (14) -- (16);
			\draw (15) -- (17);
			\draw (12) -- (14);
			\draw (3) -- (16);
			\draw (3) -- (17);
			\draw (4) -- (18);
			\draw (4) -- (17);
			\draw (18) -- (5);		
			
		\end{tikzpicture} 
		 
		\vspace{0.3cm} Figure 6. ~ $P_{3} \Box C_{5}$ \\
		\label{fig 7}	
	\end{center}
	
	In $P_{3} \Box C_{5}$, 
		\begin{align*}
		N_{1}^{(1)}~&=~N(u_{1}^{(1)})~=~N(u_1)~=~\{v_1,v_2,v_5\}, \\
		N_{2}^{(2)}~&=~N( u_{2}^{(2)})~=~N(u_2)~=~\{v_2,v_5,v_4,v_7\}, \\
		N_{3}^{(3)}~&=~N(u_{3}^{(3)})~=~N(u_3)~=~\{v_4,v_7,v_3\},\\
		N_{1}^{(3)}~&=~N(u_{1}^{(3)})~=~N(u_4)~=~\{v_5,v_8,v_7\},\\
		N_{2}^{(4)}~&=~N(u_{2}^{(4)})~=~N(u_5)~=~\{v_8,v_7,v_9,v_3\},\\
		N_{3}^{(5)}~&=~N(u_{3}^{(5)})~=~N(u_6)~=~\{v_9,v_3,v_6\}.
			\end{align*}
	\begin{align*}
	N_{1}^{(1)} \cap  N_{2}^{(2)} = \{v_2,v_5\} \neq \emptyset,~
	N_{2}^{(2)} \cap N_{3}^{(3)} = \{v_4,v_7\} \neq \emptyset.\\
	N_{1}^{(3)} \cap N_{2}^{(4)} = \{v_7,v_8\} \neq \emptyset,~
	N_{2}^{(4)} \cap N_{3}^{(5)} = \{v_3,v_9\} \neq \emptyset.
\end{align*}

$\Rightarrow$ $N_{1}^{(1)} N_{2}^{(2)} N_{3}^{(3)}$ \& $N_{1}^{(3)} N_{2}^{(4)} N_{3}^{(5)}$ are nbh chains in $P_{3} \times C_{5}$. Also,
	  \begin{align}
	 	N_{1}^{(1)} \setminus  N_{2}^{(2)} ~&=~\{v_1\} ~\& ~N_{3}^{(3)} \setminus N_{2}^{(2)} ~=~\{v_3\}\label{b_1};\\
	 		N_{1}^{(3)} \setminus  N_{2}^{(4)} ~&=~\{v_5\} ~\& ~N_{3}^{(5)} \setminus N_{2}^{(4)} ~=~\{v_6\}\label{b_2};\\
	 		N_{2}^{(2)} \setminus N_{1}^{(1)} ~&=~\{v_4,v_7\}~\subset ~ \{v_4,v_7,v_3\}~= N_{3}^{(3)};\label{c_1}\\
	 	N_{2}^{(4)} \setminus  N_{1}^{(3)} ~&=~\{v_3,v_9\}~\subset \{v_3,v_9,v_6\} = N_{3}^{(5)}; \label{d_1}\\
	 	\vert  N_{2}^{(2)} \setminus N_{1}^{(1)} \vert &~=~2 >~1~ and~ ~\vert N_{2}^{(4)} \setminus  N_{1}^{(3)} \vert~=~2 >~1.\label{e_1} 		
	 \end{align}
	 Equations (\ref{b_1}) to (\ref{e_1}) implies, $N_{1}^{(1)} N_{2}^{(2)} N_{3}^{(3)}$ \& $N_{1}^{(3)} N_{2}^{(4)} N_{3}^{(5)}$ are Type-1 nbh chains in $P_{3} \Box C_{5}$.
	 
	To complete the proof, we have to  show that the above two nbh chains satisfy conditions of NC-T2. Consider the following.
	\begin{align}
		N_{1}^{(1)} \cap N_{2}^{(2)} \cap N_{1}^{(3)} = \{v_5\} = N_{1}^{(3)} \setminus N_{2}^{(4)} \neq \emptyset; \label{f}\\
		N_{3}^{(3)} \cap N_{3}^{(5)} \cap N_{2}^{(4)} = \{v_3\} = N_{3}^{(3)} \setminus N_{2}^{(2)} \neq \emptyset; \label{g}\\	
		N_{2}^{(2)} \cap N_{3}^{(3)} \cap N_{2}^{(4)} \cap N_{1}^{(3)} = \{v_7\} \neq \emptyset.\label{j} 
	\end{align}

	
	\begin{center}
		\begin{tikzpicture}  
			[scale=.8,auto=center,every node/.style={draw,circle,fill=blue!20}] 

			\node (1) at (-8,2)[circle,fill=red!50][label=0:$\tiny u_{1}^{(1)}$]  {$u_{1}=$};
			
			\node (2) at (-5,2)[circle,fill=red!50][label=0:$\tiny u_{2}^{(2)} $] {$u_{2}=$};
			\node (3) at (-9.5,4)  [label=0:$\tiny u_{1}^{(5)}$]{$v_{1}=$};
			\node (4) at (-6.5,4)  [label=0:$\tiny u_{2}^{(1)}$]{$v_{2}=$};
			\node (5) at (-3.5,4) [label=0:$\tiny u_{3}^{(2)}$] {$v_{4}=$};
			\node (7) at (-2,2) [circle,fill=red!50] [label=0:$\tiny u_{3}^{(3)}$]{$u_{3}=$};
				
			\node (10) at (-6.5,0)[label=0:$\tiny u_{1}^{(2)} $] {$v_{5}=$};
			\node (11) at (-3.5,0) [label=0:$\tiny u_{2}^{(3)}$] {$v_{7}=$};
			\node (12) at (-0.5,0) [label=0:$\tiny u_{3}^{(4)} $]{$v_{3}=$};
			\node (13) at (-5,-2) [circle,fill=red!50] [label=0:$\tiny u_{1}^{(3)}$] {$u_{4}=$};
			\node (14) at (-2,-2) [circle,fill=red!50] [label=0:$\tiny u_{2}^{(4)}$] {$u_{5}=$};
			\node (15) at (1,-2)  [circle,fill=red!50] [label=0:$\tiny u_{3}^{(5)}$]{$u_{6}=$};
			\node (16) at (-3.5,-4) [label=0:$\tiny u_{1}^{(4)}$] {$v_{8}=$};
			\node (17) at (-0.5,-4) [label=0:$\tiny u_{2}^{(5)}$] {$v_{9}=$};
			\node (18) at (2.5,-4)  [label=0:$\tiny 
			u_{3}^{(1)}$]{$v_{6}=$};	
			
			\draw (1) -- (3);
			\draw (1) -- (4);
			\draw (1) -- (10);
			\draw (2) -- (4);
			\draw (2) -- (5);
			\draw (2) -- (10);
			\draw (2) -- (11);
			\draw (7) -- (5);
			\draw (7) -- (11);
			
			\draw (7) -- (12);
			\draw (10) -- (13);
			\draw (11) -- (13);
			\draw (11) -- (14);
			\draw (12) -- (15);
			\draw (13) -- (16);
			\draw (14) -- (17);
			\draw (15) -- (18);	
			\draw (14) -- (16);
			\draw (15) -- (17);
			\draw (12) -- (14);
			\draw[dashed]
			(3) to[out=-110,in=-150] (17);
			
			\draw[dashed]
			(3) to[out=-90,in=150] (16);
			\draw[dashed]
			(4) to[out=-73,in=145] (17);
			\draw[dashed]
			(4) to[out=-20,in=150] (18);
			\draw[dashed]
			(5) to[out=-20,in=80] (18);			
	\end{tikzpicture}  
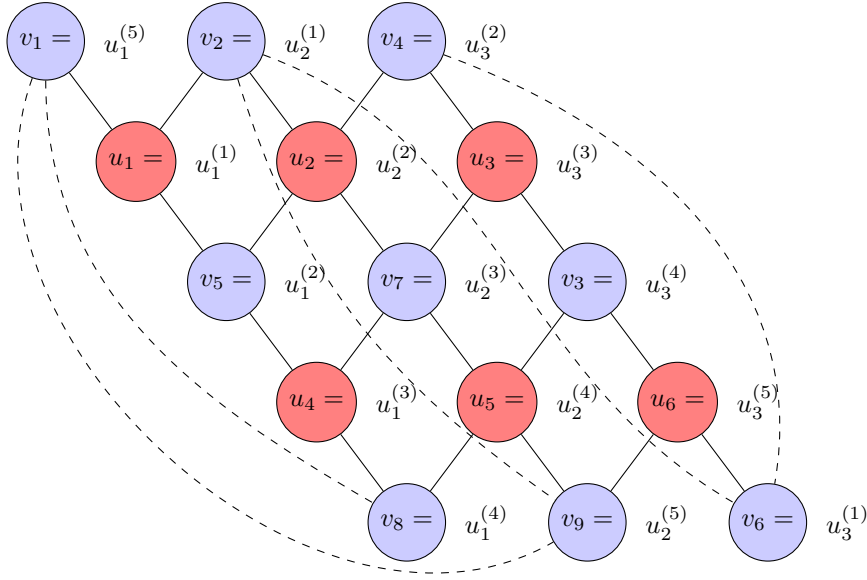
		 Figure 7. ~$P_{3} \Box C_{5}$ 
	\end{center}

	$\Rightarrow$ Nbh chains $N_{1}^{(1)} N_{2}^{(2)} N_{3}^{(3)}$ \& $N_{1}^{(3)} N_{2}^{(4)} N_{3}^{(5)}$ are Type-1 and they form a Type 2 nbh chains in $P_3 \Box C_5$.

In Figure 7, we emphasis nbh chains of Type-1   \small $N_{1}^{(1)} N_{2}^{(2)} N_{3}^{(3)}$ \& $N_{1}^{(3)} N_{2}^{(4)} N_{3}^{(5)}$\normalsize which form NC-T2 in $P_3 \Box C_5$. For more clarity, similar to Figure 5, we draw Figure 7 where the two nbh chains are $N_{1}^{(1)} N_{2}^{(2)} N_{3}^{(3)}$ \& $N_{1}^{(3)} N_{2}^{(4)} N_{3}^{(5)}$ which form NC-T2 in $P_3 \Box C_5$. See Figure 7.
\hfill $\Box$ 

\begin{theorem} \quad \label{t1} {\rm For $k \geq 2$, $n \geq 3$, graph $P_{k} \Box C_{n}$ contains NC-T2, $n,k \in \mathbb{N}$.}
\end{theorem}
\begin{proof}\quad Let $G = P_{k} \Box C_{n}$, $k \geq 2$ and $n \geq 3$. $G$ contains $n$ copies of $P_{k}$ and $k$ copies of $C_{n}$. Let  $P_{k}^{(j)}$ = $u_{1}^{(j)} u_{2}^{(j)} \ldots u_{k}^{(j)}$ be the $j^{th}$ copy of $P_{k}$ and $C_{n}^{(i)}$ = $(u_{i}^{(1)}u_{i}^{(2)}\ldots u_{i}^{(n)})$ be the $i^{th}$ copy of $C_{n}$ in $G$, $1 \leq i \leq k$ and $1 \leq j \leq n$, $ u_{i}^{(j)} \in V(G) $, $ i $ is under arithmatic modulo $ k $ and $ j $ is under arithmatic modulo $ n $.\\
	Let $N_{1}^{(j)}$ = $N(u_{1}^{(j)})$ = $\left\{u_{1}^{(j-1)},u_{1}^{(j+1)},u_{2}^{(j)}\right\}  $, $N_{k}^{(j)}$ = $N(u_{k}^{(j)})$ = $\left\{u_{k-1}^{(j)},u_{k}^{(j-1)},u_{k}^{(j+1)}\right\}$ and $N_{i}^{(j)}$ = $N(u_{i}^{(j)})$ = $\left\{u_{i-1}^{(j)},u_{i+1}^{(j)},u_{i}^{(j-1)},u_{i}^{(j+1)}\right\}  $, $1 < i < k$, $ 1 \leq j \leq n $.
	
	$\Rightarrow$ $N_{i}^{(j)} \cap N_{i+1}^{(j+1)}$ $ = \{u_{i}^{(j+1)},u_{i+1}^{(j)}\} $ $ \neq$ $\emptyset$ for $1 \leq i \leq k-1$ \& $ 1 \leq j \leq n $.
	
	$\Rightarrow$ $ N_{1}^{(j)}N_{2}^{(j+1)}\cdots N_{k}^{(j+(k-1))}$  is a nbh chain contained in $G$, $ 1 \leq j \leq n $.
	
	Also, $N_{i+1}^{(j+1)} \setminus N_{i}^{(j)} \subset N_{i+2}^{(j+2)}$ for  $1 \leq i \leq k-2$, $ 1 \leq j \leq n $ and
	
	$\vert N_{i}^{(j)} \setminus N_{i+1}^{(j+1)} \vert$, $\vert N_{i+1}^{(j+1)} \setminus N_{i}^{(j)} \vert \geq 1$ for $1 \leq i \leq k-1$, $ 1 \leq j \leq n $.
	
	Thus the nbh chain $ N_{1}^{(j)}N_{2}^{(j+1)}\cdots N_{k}^{(j+(k-1))} $ satisfies condition $(ii)$ of NC-T1 for  $1 \leq j \leq n$. 
	
	Consider  the two nbh chains  $N_{1}^{(1)} N_{2}^{(2)} \cdots N_{k}^{(k)}$ and  $N_{1}^{(3)} N_{2}^{(4)} \cdots N_{k}^{(k+2)}$ in $ G $.
	
	We have, $N_{1}^{(1)} \setminus N_{2}^{(2)}$ = $\left\{u_{1}^{(n)}\right\}$ and 	$N_{k}^{(k)}$ $\setminus$ $N_{k-1}^{(k-1)} $= $\{ u_{k}^{(k+1)} \}$;
	
	$N_{1}^{(3)} \setminus N_{2}^{(4)}$ = $\left\{u_{1}^{(2)}\right\}$ and  $N_{k}^{(k+2)}$ $\setminus$ $N_{k-1}^{(k+1)}$ = $\{ u_{k}^{(k+3)} \}$.
	
	$\Rightarrow$  $N_{1}^{(1)} N_{2}^{(2)} \cdots N_{k}^{(k)}$ and $N_{1}^{(3)} N_{2}^{(4)} \cdots N_{k}^{(k+2)}$ satisfy condition $(i)$ of NC-T1.
	
	$\Rightarrow$  $N_{1}^{(1)} N_{2}^{(2)} \cdots N_{k}^{(k)}$ and  $N_{1}^{(3)} N_{2}^{(4)} \cdots N_{k}^{(k+2)}$  are Type-1 nbh chains in $ G$.
	
	To complete the proof, we have to  show that the above two nbh chains satisfy conditions of NC-T2. Consider the following.
	\begin{enumerate} 
		\item [\rm (i)]  $N_{1}^{(1)} \cap N_{2}^{(2)} \cap N_{1}^{(3)} = \left\{u_{1}^{(2)}\right\}$ \& $N_{1}^{(3)}\setminus N_{2}^{(4)} = \left\{u_{1}^{(2)}\right\}$.
		
		$\Rightarrow$ $N_{1}^{(1)} \cap N_{2}^{(2)} \cap N_{1}^{(3)}$ = $N_{1}^{(3)}\setminus N_{2}^{(4)} \neq \emptyset$.  
		
		Also,	$N_{k}^{(k)} \cap N_{k}^{(k+2)}  \cap N_{k-1}^{(k+1)}$ = $\left\{u_{k}^{(k+1)}\right\}$ \& $N_{k}^{(k)} \setminus N_{k-1}^{(k-1)}$ = $\left\{u_{k}^{(k+1)}\right\}$. 
		
		$\Rightarrow$	$N_{k}^{(k)} \cap N_{k}^{(k+2)}  \cap N_{k-1}^{(k+1)}$= $N_{k}^{(k)} \setminus N_{k-1}^{(k-1)} \neq \emptyset$.
		
		\item [\rm (ii)]  $N_{i}^{(i)} \cap N_{i+1}^{(i+1)} \cap N_{i}^{(i+2)} \cap N_{i-1}^{(i+1)} = \left\{u_{i}^{(i-1)}, u_{i}^{(i+1)}, u_{i-1}^{(i)}, u_{i+1}^{(i)}\right\}$ 
		
		\hfill	$\cap \left\{u_{i+1}^{(i)}, u_{i+1}^{(i+2)}, u_{i}^{(i+1)}, u_{i+2}^{(i+1)}\right\}$ 	$\cap \left\{u_{i}^{(i+1)}, u_{i}^{(i+3)}, u_{i-1}^{(i+2)}, u_{i+1}^{(i+2)}\right\}$ 
		
		$\cap \left\{u_{i-1}^{(i)}, u_{i-1}^{(i+2)}, u_{i-2}^{(i+1)}, u_{i}^{(i+1)}\right\}$  = $\left\{u_{i}^{(i+1)}\right\} \neq \emptyset$, $1 < i < k$.
	\end{enumerate}
	Thus, the two nbh chains  $N_{1}^{(1)} N_{2}^{(2)} \cdots N_{k}^{(k)}$ and  $N_{1}^{(3)} N_{2}^{(4)} \cdots N_{k}^{(k+2)}$ form a Type 2 nbh chains in $G$. Hence the result.
\end{proof}

\begin{cor} \quad \label{tt1} {\rm For $k \geq 2$ and $n \geq 3$, graph $P_{k} \Box C_{n}$ contains two NC-T1, each of length $k$ and the two form NC-T2, $n,k \in \mathbb{N}$. \hfill $\Box$}
\end{cor}

	\vspace{.2cm}	
	The following theorem establishes a certain family of graphs as NDM by the presence of NC-T1 in these graphs.
	
	\begin{theorem}\quad \label{t2} {\rm Let $G$ be a graph containing Type-1 nbh chain of length $2n$, $n\in\mathbb{N}$. Then $G$ is NDM.	}
	\end{theorem} 	
	\begin{proof} Let $N_1 N_2 \cdots N_{2n}$ be a nbh chain of Type-1 in $G$. We prove the theorem by the method of contradiction.
		
		If possible, let $G$ be a DM graph. Let $f$ be a DML of $G$ with DM constant $S$. For sets $N_1$ and $N_2$, we get, 
			\begin{align}
			~~	N_{1}~&=~ (N_{1} \setminus N_2) \cup (N_{1}\cap N_{2}), \label{1}\\ ~~~	N_{2}~&= (N_{1}\cap N_2) \cup ~ (N_{2} \setminus N_1).\label{2}
		\end{align}
		Using the definition of DML, we get,
		\begin{align}
			S~=\sum_{v_j \in  N_1}f(v_j)=~\sum_{v_j \in N_2 }f(v_j)~=~\cdots~=~\sum_{v_j \in  N_i}f(v_j)\notag \\ = \cdots~=~\sum_{v_j \in N_{2n-1}}f(v_j)~~=~~\sum_{v_j \in N_{2n} }f(v_j). \label{3}	\end{align}
		
		$$ (\ref{1}),(\ref{2})~\&~ (\ref{3}) \implies 
		\sum_{v_j \in (N_1 \setminus N_2)~\cup~(N_1 \cap N_2)}f(v_j) = \sum_{v_j \in (N_1 \cap N_2)~\cup~(N_2 \setminus N_1)}f(v_j). $$
		
		\begin{align}
			\Rightarrow \sum_{v_j \in N_1 \setminus N_2}f(v_j)+ \sum_{v_j \in N_1 \cap N_2}f(v_j)~~=~~\sum_{v_j \in N_1 \cap N_2}f(v_j) + \sum_{v_j \in N_2 \setminus N_1}f(v_j). \notag
		\end{align}
		\begin{align}
			\Rightarrow \sum_{v_j \in N_1 \setminus N_2}f(v_j)~&=~~\sum_{v_j \in N_2 \setminus N_1}f(v_j). \notag 
		\end{align}
		$$ \Rightarrow  
		\sum_{v_j \in N_1 \setminus N_2}f(v_j)~~=f(v_1)~=~~\sum_{v_j \in N_2 \setminus N_1}f(v_j) ~ \text{using property of NC-T1}. \label{5} $$
		\begin{center}
			\begin{tikzpicture}[thick,
				set/.style = { circle, minimum size = .061cm}]
				
				\draw[green] (3,0) ellipse (2cm and 1cm);
				\draw[blue] (1,0) ellipse (1.94cm and 0.99cm);
				
				\begin{scope}[red]
					\clip(2,-1)rectangle(3,1.2);
					\draw (1,0) circle [x radius=2cm, y radius=10mm];
				\end{scope}
				\begin{scope}[red]
					\clip(2,-1)rectangle(3,1.2);
					\draw (1,0) circle [x radius=1.975cm, y radius=10mm];
				\end{scope}
				\begin{scope}[blue]
					\clip(2,-1)rectangle(3,1.2);
					\draw (1,0) circle [x radius=1.935cm, y radius=9.8mm];
				\end{scope}

				\begin{scope}[red]
					\clip(2,-1.5)rectangle(7,1.5);
					\draw (4,0) circle [x radius=3cm, y radius=11.9mm];
				\end{scope}
				
				\node[blue] at (0,0) {$N_{i}\setminus N_{i+1}$};
				\node at (2,0) {\tiny $N_{i}\cap N_{i+1}$};
				\node[green] at (4,0.2) {$N_{i+1}\setminus N_{i}$};
				\node[red] at (6,-0.2) {$N_{i+2}\setminus N_{i+1}$};	
			\end{tikzpicture}\\
			\vspace{0.2cm} Figure $8$ \label{fig 9}			
		\end{center}
	
		In Figure 8, closed curves with colors blue, green and red represent nbhs $N_i$, $N_{i+1}$, $N_{i+2}$, respectively. And property (ii) of NC-T1, $ (N_{i+1} \setminus N_{i}) \subset N_{i+2}$, is presented in the figure. See Figure 8.
		
		For $i = 1,2,\ldots,2n-2$,  
		\begin{align}
			(N_{i+1} \setminus N_{i}) \subset N_{i+2}.	\Rightarrow ~N_{i+1} \setminus N_{i}~ = N_{i+1} \cap N_{i+2}.\label{6}
		\end{align}
		
		Also, for $ i ~=~1,2,\ldots,2n-1 $, we have,
		\begin{align}
			N_{i}&= (N_{i} \setminus N_{i+1}) \cup (N_{i}\cap N_{i+1})~ \text{and}\label{7} \\
			N_{i+1}&= (N_{i+1} \setminus N_{i}) \cup (N_{i}\cap N_{i+1}).\label{8} 
		\end{align}
		
		This implies, for  $ i ~=~1,2,\ldots,2n-2$, 
		\begin{align}
			N_{i+2}&= (N_{i+2} \setminus N_{i+1}) \cup (N_{i+1}\cap N_{i+2})\notag\\
			&= (N_{i+2} \setminus N_{i+1}) \cup (N_{i+1} \setminus N_{i}) ~\text{using (\ref{6})}.\label{9}
		\end{align}
		\begin{align} 
			\Rightarrow		N_{i+3}&= (N_{i+3}\setminus N_{i+2}) \cup (N_{i+2} \setminus N_{i+1}) ~ \text{for}~ i= 1,2,\dots ,2n-3. \label{10} 
		\end{align}
		Using (\ref{7})~\&~(\ref{8})~ in  (\ref{3}), for  $ i=1,2,\dots ,2n-1, $ ~we get,
		\begin{align*}	S = \sum_{u \in N_{i} \setminus N_{i+1}}{f(u)}+\sum_{u\in N_{i}\cap N_{i+1}}{f(u)}=\sum_{u\in N_{i+1} \setminus N_{i}}{f(u)}+\sum_{u\in N_{i}\cap N_{i+1}}{f(u)}.
		\end{align*}
		\begin{align}
			\Rightarrow	 \sum_{u\in N_{i} \setminus N_{i+1}}{f(u)}~~=\sum_{u\in N_{i+1} \setminus N_{i}}{f(u)}~ for ~ i = 1,2,\dots ,2n-1. \label{11}
		\end{align}
		Similarly, using (\ref{9})~\&~(\ref{10})~ in  (\ref{3}), for $ i=1,2,\dots ,2n-3,$ ~we get,
		$$	S = \sum_{u\in N_{i+2} \setminus N_{i+1}}{f(u)} + \sum_{u\in N_{i+1} \setminus N_{i}}{f(u)} \hspace{4cm}	$$
		$$	\hspace{3cm}	= \sum_{u\in N_{i+3} \setminus N_{i+2}}{f(u)}~+\sum_{u\in N_{i+2} \setminus N_{i+1}}{f(u)}.$$
		\begin{align}
			\Rightarrow	\sum_{u\in N_{i+1} \setminus N_{i}}{f(u)} ~ 
			=\sum_{u\in N_{i+3} \setminus N_{i+2}}{f(u)}~ \text{for} ~ i = 1,2,\dots ,2n-3.\label{120}
		\end{align}	
		\begin{align}
			\Rightarrow \sum_{u\in N_{2} \setminus N_{1}}{f(u)} ~ 
			= \sum_{u\in N_{4} \setminus N_{3}}{f(u)} = . . .  = \sum_{u\in N_{2n} \setminus N_{2n-1}}{f(u)}. \label{12}
		\end{align}
		Using the definition of NC-T1 in (\ref{11}), we get,  	
		\begin{align*}
			\sum_{u\in N_{1} \setminus N_{2}}{f(u)} = \sum_{u\in N_{2} \setminus N_{1}}{f(u)} = f(v_1)~ and \hspace{3.5cm} (a) 
		\\	
			\sum_{u\in N_{2n-1} \setminus N_{2n}}{f(u)} = \sum_{u\in N_{2n} \setminus N_{2n-1}}{f(u)} = f(v_{2n}). \hspace{3.5cm} (b)
		\end{align*}
	
		$$(\ref{12}), (a), (b)	\Rightarrow f(v_1) = 	\sum_{u\in N_{1} \setminus N_{2}}{f(u)} = \sum_{u\in N_{2} \setminus N_{1}}{f(u)}  = \sum_{u\in N_{4} \setminus N_{3}}{f(u)}$$
		$$ = . . .  = \sum_{u\in N_{2n} \setminus N_{2n-1}}{f(u)} = f(v_{2n}).$$
			$\Rightarrow f(v_1) = 	f(v_{2n})$ which  is a contradiction to DML.
			
		Hence the result. 
	\end{proof}
	
	\begin{cor} \quad \label{c1}  {\rm For $n \geq 3$ and $n,k \in \mathbb{N}$, graph $P_{2k} \Box C_{n}$ is NDM.}	
	\end{cor}
	\begin{proof} \quad  Proof follows from Corollary \ref{tt1} and Theorem \ref{t2}.  
	\end{proof}
Thus, we are left with the following conjecture instead of Conjecture \ref{ca1}.
	
	\begin{conj} \quad \label{c2} {\rm For $n \geq 3$ and $n,k \in \mathbb{N}$, graph $P_{2k+1} \Box C_{n}$ is NDM. \hfill $\Box$	}
	\end{conj}
	
	\section{Families of NDM graphs by existence of NC-T1}
	
 Using Theorem \ref{t2}, we present some families of graphs which are NDM. At first, let us start with a few definitions which are required in this section.
		
	\begin{definition}\cite{xw} For $k\in\mathbb{N}$ and $k \geq 2$,
		let $\{G_i\}^k_{i=1}$ be a collection of graphs with $u_i\in V(G_i)$ as a fixed vertex, $1 \leq i \leq k$. \em{The vertex amalgamation}, denoted
		by $Amal(G_i, \{u_i\}, k)$, is a graph formed by taking all the $G_i$'s and identifying $u_i$'s. 
	\end{definition}
	
	\begin{definition}\cite{ep92} A \emph{two terminal graph} is a graph with two distinguished vertices, $s$ and $t$ called \emph{source} and \emph{sink}, respectively. 
	\end{definition}
	
	\begin{definition}\cite{ep92} \emph{The series composition of 2 two terminal graphs} $X$ and $Y$ is a two terminal graph created from the disjoint union of graphs $X$ and $Y$ by merging the sink of $X$ with the source of $Y$. The source of $X$  becomes the source of series composition and the sink of $Y$ becomes the sink of series composition. 
	\end{definition}
	
	\begin{definition}\cite{ep92} \emph{The parallel composition of two terminal graphs} $X$ and $Y$ is a two terminal graph created from the disjoint union of graphs $X$ and $Y$ by merging the sources of $X$ and $Y$ to create the source of parallel composition and merging the sinks of $X$ and $Y$ to create the sink of parallel composition. 
	\end{definition}
	
	\begin{definition}\cite{go78} Let $G$ be a graph of order $ n(G)$ and $H$ be a graph with root vertex $v$. Then, the \emph{rooted product graph} of $G$ and $H$ is defined as the graph obtained from $G$ and $H$ by taking one copy of $G$ and $n(G)$ copies of $H$ and identifying the $i^{th}$ vertex of $G$  with the root vertex $v$ in the $i^{th}$ copy of $H$ for every $i\in \{ 1,2,3,...,n(G) \}$.
	\end{definition}

The following families of graphs are NDM by the existence of NC-T1 of even length in each of these graphs.
	
  \begin{enumerate}
		
   \item	Let $ G_1,G_2,\cdots,G_k $ be mutually disjoint graphs such that at least one of them contain a NC-T1 of even length. Then $\displaystyle$ $ \bigcup_{i=1}^k G_i $ is NDM.
		
   \item Let $k\in\mathbb{N}$, $k \geq 2$ and $\{G_i\}^k_{i=1}$ be a collection of mutually disjoint graphs with $u_i\in V(G_i)$ as a fixed vertex such that at least one of them contain an NC-T1 of even length. Let $H$ = $Amal(G_i, \{u_i\}, k)$, the vertex amalgamation of $ G_1,G_2,\cdots,G_k$ by taking all the $G_i$'s and identifying $u_i$'s without disturbing the structure of at least one NC-T1 of even length in $H$, $1 \leq i \leq k$. Then $H$ is NDM.
		
   \item Let $G$ be a simple graph having a NC-T1 of even length, say $NSG_1$ and $H$ be a super graph of $G$ without disturbing the structure of $NSG_1$ in $H$ and $NSG_1$ is also an NC-T1 in $H$. Then $H$ is NDM.
		
  \item	Let $ G $ be a simple graph containing an NC-T1 of even length, say $NSG_1$. Then any subgraph $H$ of $G$ containing $NSG_1$ as a NC-T1 of even length is also NDM.
		
  \item	Let $G$ be a two terminal simple graph and $H$ be a graph having NC-T1 of even length. Then the series composition of $G$ and $H$ is NDM, follows from the structure of the combined graph.
		
  \item	Let $G$ be a two terminal simple graph and $H$ be a graph having NC-T1 of even length. Then the parallel composition of $G$ and $H$ is NDM, follows from the structure of the combined graph.
		
 \item	Let $G$ be any simple graph and $H$ be a graph having NC-T1 of even length. Then the rooted product of $G$ and $H$ is NDM.
		
\end{enumerate}
	
	\noindent
	{\bf Conclusion.} \quad In this work, we obtained families of non-distance magic graphs by the existence of NC-T1 of even length in these graphs and partially settled the conjecture on NDM on cylindrical graphs $P_{2m} \Box C_n$ for $n \geq 3$ and $m,n \in \mathbb{N}$. Recently, the authors completed the other part of the conjecture and submitted the work for publication. And the authors feel that a lot of scope is there to utilise the concepts of nbh chains and nbh sequence graphs in graphs when their order is large while dealing with graphs corresponding to large data and their network problems. Also, nbh sequence, nbh walk, nbh trail, nbh path/chain concepts can be used in the theory of hypergraphs. 
	
	\vspace{.5cm}
	\noindent
	\textbf{Conflict of Interest}
	
	\noindent
	\textit {The authors declare that there is no conflict of interests regarding the publication of this paper.}
	
	\vspace{.5cm}
	\noindent
	\textbf{Acknowledgement}
	
	\noindent
	\textit {\rm We express our sincere thanks to the Central University of Kerala, Kasaragod - 671 316, Kerala, India for providing facilities to carry out this research work.}

\end{document}